\documentclass[10pt,a4paper]{article}
\usepackage[latin1]{inputenc}
\usepackage[T1]{fontenc}
\usepackage{amsmath,amscd}
\usepackage{graphicx}
\usepackage[all]{xy}
\usepackage{amsthm}
\usepackage{amssymb,url}
\usepackage{bbm}
\usepackage{ae}
\usepackage[all]{xy}

\newtheorem{sats}{Theorem}[section]
\newtheorem{sats*}{Theorem}

\newtheorem{lem}[sats]{Lemma}

\newtheorem{prop}[sats]{Proposition}

\newcommand{\R}{\mathbbm{R}}
\newcommand{\C}{\mathbbm{C}}

\newcommand{\Z}{\mathbbm{Z}}

\newcommand{\ellL}{\mathcal{L}}

\newcommand{\rd}{\mathrm{d}}

\newcommand{\ch}{ \mathrm{cs}}
\newcommand{\che}{ \mathrm{ch}}

\newcommand{\Ko}{\mathcal{K}}
\newcommand{\Bo}{\mathcal{B}}

\newcommand{\He}{\mathcal{H}}

\newcommand{\ind}{\mathrm{i} \mathrm{n} \mathrm{d} \,}

\renewcommand{\epsilon}{\varepsilon}
\renewcommand{\phi}{\varphi}
\renewcommand{\Re}{\mathrm{Re}\,}

\newcommand{\e}{\mathrm{e}}
\newcommand{\tra}{\mathrm{t}\mathrm{r}}

\title{Analytic formulas for topological degree of non-smooth mappings: the odd-dimensional case}
\author{Magnus Goffeng}
\date{Department of Mathematical Sciences, Division of Mathematics\\
Chalmers university of Technology and University of Gothenburg}

\begin{document}
\maketitle

\begin{abstract}
The notion of topological degree is studied for mappings from the boundary of a relatively compact strictly pseudo-convex domain in a Stein manifold into a manifold in terms of index theory of Toeplitz operators on the Hardy space. The index formalism of non-commutative geometry is used to derive analytic integral formulas for the index of a Toeplitz operator with H\"older continuous symbol. The index formula gives an analytic formula for the degree of a H\"older continuous mapping from the boundary of a strictly pseudo-convex domain.
\end{abstract}

\section*{Introduction}

This paper is a study of analytic formulas for the degree of a mapping from the boundary of a relatively compact strictly pseudo-convex domain in a Stein manifold. The degree of a continuous mapping between two compact, connected, oriented manifolds of the same dimension is abstractly defined in terms of homology for continuous functions. If the function $f$ is differentiable, an analytic formula can be derived using Brouwer degree, see \cite{reyeche}, or the more global picture of de Rham-cohomology. For any form $\omega$ of top degree the form $f^*\omega$ satisfies 
\[\int _X f^*\omega=\deg f \int _Y\omega.\] 
Without differentiability conditions on $f$, there are no known analytic formulas beyond the special case of a H\"older continuous mapping $S^1\to S^1$ which can be found in Chapter $2.\alpha$ of \cite{connes}. The degree of a H\"older continuous function $f:S^1\to S^1$ of exponent $\alpha$ is expressed by an analytic formula by replacing de Rham cohomology with the cyclic homology of the algebra of H\"older continuous functions as 
\begin{equation}
\label{connesalpha}
\deg(f)=\frac{1}{(2\pi i)^{2k}}\int f(z_0)\frac{f(z_1)-f(z_0)}{z_1-z_0}\cdots \frac{f(z_0)-f(z_{2k})}{z_0-z_{2k}}\rd z_0\ldots \rd z_{2k},
\end{equation}
whenever $\alpha(2k+1)>1$. Later, the same technique  was used in \cite{grigone} and \cite{grigtwo} in constructing index formulas for pseudo differential operators with operator valued symbols. Our aim is to find new formulas for the degree in the multidimensional setting by expressing the degree of a H\"older continuous function as the index of a Toeplitz operator and using the approach of \cite{connes}. 

The motivation to calculate the degree of a non-smooth mapping comes from non-linear $\sigma$-models in physics. For instance, the Skyrme model describing self-interacting mesons in terms of a field $f:X\to Y$, see \cite{auckkapi}, only have a constant solution if one does not pose a topological restriction and since the solutions are rarely smooth, but rather in the Sobolev space $W^{1,d}(X,Y)$, one needs a degree defined on non-continuous functions. In the paper \cite{brenir}, the notion of a degree was extended as far as to VMO-mappings in terms of approximation by continuous mappings. See also \cite{breli} for a study of the homotopy structure of $W^{1,d}(X,Y)$. 

The main idea that will be used in this paper is that the cohomological information of a continuous mapping $f:X\to Y$ between odd dimensional manifolds can be found in the induced mapping $f^*:K^1(X)\to K^1(Y)$ using the Chern-Simons character. The analytic formula will be obtained by using index theory of Toeplitz operators. The index theory of Toeplitz operators is a well studied subject for many classes of symbols, see for instance \cite{badog}, \cite{boutetdemonvel}, \cite{connes} and \cite{guehig}. If $X=\partial\Omega$, where $\Omega$ is a strictly pseudo-convex domain in a complex manifold, and $f:\partial\Omega\to Y$ is a smooth mapping the idea can be expressed by the commutative diagram:
\begin{equation}
\label{smodia}
\begin{CD}
K^1(Y) @>f^*>>K^1(\partial\Omega)@>\ind >> \Z \\
@VV\ch_YV @VV\ch_{\partial\Omega}V @VVV\\
H^{odd}_{dR}(Y)@>f^*>>H^{odd}_{dR}(\partial\Omega) @>\chi_{\partial\Omega} >> \C \\
\end{CD} 
\end{equation}
where the mapping $\ind:K^1(\partial\Omega)\to \Z$ denotes the index mapping defined in terms of suitable Toeplitz operators on $\partial\Omega$ and 
\[\chi_{\partial\Omega}(x):=-\int_{\partial\Omega} x\wedge Td(\Omega).\] 
The left part of the diagram \eqref{smodia} is commutative by naturality of the Chern-Simons character and the right part of the diagram is commutative by the Boutet de Monvel index formula. 

The $K$-theory is a topological invariant and the picture of the index map as a pairing in a local homology theory via Chern-Simons characters can be applied to more general classes of functions than the smooth functions. The homology theory present through out all the index theory is cyclic homology. For a H\"older continuous mapping $f:\partial\Omega\to Y$ of exponent $\alpha$ and $\Omega$ being a relatively compact strictly pseudo-convex domain in a Stein manifold the analogy of the diagram \eqref{smodia} is
\begin{equation}
\label{holdia}
\begin{CD}
K_1(C^\infty(Y))@>f^*>>K_1(C^\alpha(\partial\Omega))@>\ind >> \Z \\
@VV\ch_YV @VV\ch_{\partial\Omega}V @VVV\\
HC_{odd}(C^\infty(Y))@>f^*>>HC_{odd}(C^\alpha(\partial\Omega)) @>\tilde{\chi}_{\partial\Omega} >> \C \\
\end{CD} 
\end{equation}
where the mapping $\tilde{\chi}_{\partial\Omega}:HC_{odd}(C^\alpha(\partial\Omega))\to \C$ is a cyclic cocycle on $C^\alpha(\partial\Omega)$ defined as the Connes-Chern character of the Toeplitz operators on $\partial\Omega$, see more in \cite{connes} and \cite{connesncdg}. The condition on $\Omega$ to lie in a Stein manifold ensures that the cyclic cocycle $\tilde{\chi}_{\partial\Omega}$ can be defined on H\"older continuous functions, see below in Theorem \ref{cycref}. The right hand side of the diagram \eqref{holdia} is commutative by Connes' index formula, see Proposition $4$ of Chapter IV$.1$ of \cite{connes}. The dimension in which the Chern-Simons character will take values depends on the H\"older exponent $\alpha$. More explicitly, the cocycle $\tilde{\chi}_{\partial\Omega}$ can be chosen as a cyclic $2k+1$-cocycle for any $2k+1>2n/\alpha$.

The index of a Toeplitz operator $T_u$ on the vector valued Hardy space $H^2(\partial\Omega)\otimes \C^N$ with smooth symbol $u:\partial\Omega\to \mathrm{GL}_N(\C)$  can be calculated using the Boutet de Monvel index formula as $\ind T_u=-\int_{\partial\Omega}\ch_{\partial\Omega} [u]$ if the Chern-Simons character $\ch_{\partial\Omega} [u]$ only contains a top degree term. In particular, if $g:Y\to \mathrm{GL}_N(\C)$ satisfies that all terms, except for the top-degree term, in $\ch_{\partial\Omega} [g]$  are exact and $f:\partial \Omega\to Y$ is smooth we can consider the matrix symbol $g\circ f$ on $\partial \Omega$. Naturality of the Chern-Simons character implies the identity 
\[\deg f\int_Y\ch_{Y}[g]=-\ind T_{g\circ f}\] 
where $T_{g\circ f}$ is a Toeplitz operator on $H^2(\partial\Omega)\otimes \C^N$ with symbol $g\circ f$. This result extends to H\"older continuous functions in the sense that if we choose $g$ which also satisfies the condition $\int_Y\ch_{Y}[g]=1$ we obtain the analytic degree formula: 
\[\deg f=\tilde{\chi}_{\partial\Omega}(\ch_{\partial\Omega}[g\circ f]).\]

A drawback of our approach is that it only applies to boundaries of strictly pseudo-convex domains in Stein manifolds. We discuss this drawback at the end of the fourth, and final, section of this paper. The author intends to return to this question in a future paper and address the problem for even-dimensional manifolds. \\

The paper is organized as follows; in the first section we reformulate the degree as an index calculation using the Chern-Simons character from odd $K$-theory to de Rham cohomology. This result is not remarkable in itself, since the Chern-Simons character is an isomorphism after tensoring with the complex numbers. However, the constructions are explicit and allows us to obtain explicit expressions for a generator of the de Rham cohomology. We will use the complex spin representation of $\R^{2n}$ to construct a smooth function $u:S^{2n-1}\to SU(2^{n-1})$ such that the Chern-Simons character of $u$ is a multiple of the volume element on $S^{2n-1}$. The function $u$ will then be used to construct a smooth mapping $\tilde{g}:Y\to GL_{2^{n-1}}(\C)$ for arbitrary odd-dimensional manifold $Y$ whose Chern-Simons character coincide with $(-1)^{n}\rd V_Y$ where $\rd V_Y$ is a normalized volume form on $Y$, see Theorem \ref{degsats}. Thus we obtain for any continuous function $f:\partial\Omega\to Y$ the formula $\deg f=(-1)^{n+1}\ind T_{g\circ f}$, as is proved in Theorem \ref{indexdegreeform}.

In the second section we will review the theory of Toeplitz operators on the boundary of a strictly pseudo-convex domain. The material in this section is based on \cite{boutetdemonvel}, \cite{connes}, \cite{forn}, \cite{guehig}, \cite{henkleit} and \cite{range}. We will recall the basics from \cite{forn}, \cite{henkleit} and \cite{range} of integral representations of holomorphic functions on Stein manifolds and the non-orthogonal Henkin-Ramirez projection. We will continue the section by recalling some known results about index formulas and how a certain Schatten class condition can be used to obtain index formulas. The focus will be on the index formula of Connes, see Proposition $4$ in Chapter IV$.1$ of \cite{connes}, involving cyclic cohomology and how the periodicity operator $S$ in cyclic cohomology can be used to extend cyclic cocycles to larger algebras. In our case the periodicity operator is used to extend a cyclic cocycle on the algebra $C^\infty(\partial\Omega)$ to a cyclic cocycle on $C^\alpha(\partial \Omega)$. We will also review a theorem of Russo, see \cite{russo}, which gives a sufficient condition for an integral operator to be of Schatten class. 

The third section is devoted to proving that the Szeg\"o projection $P_{\partial\Omega}:L^2(\partial\Omega)\to H^2(\partial\Omega)$ satisfies the property that for any $p>2n/\alpha$ the commutator $[P_{\partial\Omega},a]$ is a Schatten class operator of order $p$ for any H\"older continuous functions $a$ on $\partial \Omega$ of exponent $\alpha$. The statement about the commutator $[P_{\partial\Omega},a]$ can be reformulated as the corresponding big Hankel operator with symbol $a$ being of Schatten class. We will in fact not look at the Szeg\"o projection, but rather at the non-orthogonal Henkin-Ramirez projection $P_{HR}$ mentioned above. The projection $P_{HR}$ has a particular behavior making the estimates easier and an application of Russo's Theorem implies that $P_{HR}-P_{\partial\Omega}$ is Schatten class of order $p>2n$, see Lemma \ref{henkbergdiff}. 
 
In the fourth section we will present the index formula and the degree formula for H\"older continuous functions. Thus if we let $C_{\partial\Omega}$ denote the Szeg\"o kernel and $\rd V$ the volume form on $\partial\Omega$ we obtain the following index formula for $u$ invertible and H\"older continuous on $\partial\Omega$:
\[\ind T_u=-\int_{\partial\Omega^{2k+1}} \tra\left(\prod_{i=0}^{2k}(1-u(z_i)^{-1}u(z_{i+1}))C_{\partial\Omega}(z_i,z_{i+1})\right)\rd V \]
for any $2k+1>2n/\alpha$. Here we identify $z_{2k+1}$ with $z_0$. Using the index formula for mapping degree we finally obtain the following analytic formula for the degree of a H\"older continuous mapping from $\partial\Omega$ to a connected, compact, orientable, Riemannian manifold $Y$. If $f:\partial \Omega\to Y$ is a H\"older continuous function of exponent $\alpha$, the degree of $f$ can be calculated for $2k+1>2n/\alpha$ from the formula:
\[\deg(f)=(-1)^{n}\int_{\partial\Omega^{2k+1}} \tilde{f}(z_0,z_1,\ldots, z_{2k}) \prod _{j=0}^{2k}C_{\partial\Omega}(z_{j-1},z_j)\rd V\]
where $\tilde{f}:\partial\Omega^{2k+1}\to \C$ is a function explicitly expressed from $f$, see more in equation \eqref{fsig}.

\section{The volume form as a Chern-Simons character}

In order to represent the mapping degree as an index we look for a matrix symbol whose Chern-Simons character is cohomologous to the volume form $\rd V_Y$ on $Y$. We will start by considering the case of a $2n-1$-dimensional sphere and construct a map into the Lie group $SU(2^{n-1})$ using the complex spinor representation of $Spin(\R^{2n})$. In the complex spin representation a vector in $S^{2n-1}$ defines a unitary matrix, this construction produces a matrix symbol on odd-dimensional spheres such that its Chern-Simons character spans $H^{2n-1}_{dR}(S^{2n-1})$. The matrix symbol on $S^{2n-1}$ generalizes to an arbitrary connected, compact, oriented manifold $Y$ of dimension $2n-1$ such that its Chern-Simons character coincides with $(-1)^n\rd V_Y$.

Let $V$ denote a real vector space of dimension $2n$ with a non-degenerate inner product $g$. We take a complex structure $J$ on $V$ which is compatible with the metric and extend the mapping $J$ to a complex linear mapping on $V_\C:=V\otimes_\R \C$. Since $J^2=-1$ we can decompose $V_\C:=V^{1,0}\oplus V^{0,1}$ into two eigenspaces of $J$ corresponding to the eigenvalues $\pm i$. 
If we extend $g$ to a complex bilinear form $g_\C$ on $V_\C$ and using the isomorphism $Cl_\C(V,g)\cong Cl(V_\C,g_\C)$, we can identify the complexified Clifford algebra of $V$ with the complex algebra generated by $2n$ symbols $e_{1,+},\ldots,e_{n,+},e_{1,-},\ldots,e_{n,-}$ satisfying the relations 
\[\{e_{j,+},e_{k,+}\}=\{e_{j,-},e_{k,-}\}=0\quad\mbox{and}\quad \{e_{j,+},e_{k,-}\}=-2\delta_{jk},\]
where $\{\cdot,\cdot\}$ denotes anti-commutator. The complex algebra $Cl_\C(V,g)$ becomes a $*$-algebra in the $*$-operation $e_{j,+}^*:=-e_{j,-}$.

The space $S_V:=\wedge^* V^{1,0}$ becomes a complex Hilbert space equipped with the sesquilinear form induced from $g$ and $J$. The vector space $S_V$ will be given the orientation from the lexicographic order on the basis $e_{i_1}\wedge e_{i_2}\wedge \ldots \wedge e_{i_k}$ for $i_1<i_2<\ldots < i_k$. Define $c:V_\C\to End(S_V)$ by 
\[c(v).w:=\sqrt{2}v\wedge w, \quad \mbox{for}\quad v\in V^{1,0}\quad\mbox{and}\] 
\[c(v').w:=-\sqrt{2}v'\neg w \quad \mbox{for}\quad v'\in V^{0,1}.\]
The linear mapping $c$ satisfies 
\[c(v^*)=c(v)^*\quad\mbox{and}\quad c(w)c(v)+c(v)c(w)=-2g(w,v)\] 
so by the universal property of the Clifford algebra $Cl_\C(V,g)$ we can extend $c$ to a $*$-representation $\phi:Cl_\C(V)\to End_\C(S_V)$. The space $S_V$ is a $2^{n}$-dimensional Hilbert space which we equip with a $\Z_2$-grading as follows
\[S_V=S_V^+\oplus S_V^-:=\wedge^{even}V^{1,0}\oplus \wedge^{odd} V^{1,0}.\]
Consider the subalgebra $Cl_\C(V)_+$ consisting of an even number of generators. The representation $\phi$ restricts to a representation $Cl_\C(V)_+\to End_\C(S_V^+)$ and $Cl_\C(V)_+\to End_\C(S_V^-)$. We define the $2^{n-1}$-dimensional oriented Hilbert space $E_n:=S_{\C^n}^+$ when $n$ is even and $E_n:=S_{\C^n}^-$ when $n$ is odd.  The representation $Cl_\C(\C^n)_+\to End_\C(E_n)$ will be denoted by $\phi_+$. For a vector $v\in \C^n$ we can use the fact that $\C^n\otimes_\R\C \cong \C^n\oplus \C^n$ and define 
\[v_+:=\phi_+(v\oplus 0)\in End_\C(E_n)\quad\mbox{and}\quad v_-:=\phi_+(0\oplus v)\in End_\C(E_n).\]

We will now define a symbol calculus for $S^{2n-1}$. We choose the standard embedding $S^{2n-1}\subseteq \C^n$ by taking coordinates $z_i:S^{2n-1}\to \C$ satisfying $|z_1|^2+|z_2|^2\cdots + |z_n|^2=1$. Define the smooth map $u:S^{2n-1}\to Cl_\C(\R^{2n})_+$ by 
\begin{equation}
\label{udef}
u(z):=\frac{1}{2}(e_{1,+}+e_{1,-})(z_++\bar{z}_-).
\end{equation}

\begin{prop}
\label{unitlem}
The mapping $u$ satisfies 
\[u(z)^*u(z)=u(z)u(z)^*=1\]
so $u:S^{2n-1}\to SU(2^{n-1})\subseteq End_\C(E_n)$ is well defined.
\end{prop}

The proof of this proposition is a straightforward calculation using the relations in the Clifford algebra $Cl_\C(V,g)$. Observe that if $n=2$ the mapping $u$ is a diffeomorphism since we can choose $1$ and $e_1\wedge e_2$ as a basis for $S_V^+$ and in this basis 
\[u(z_1,z_2)=\begin{pmatrix} -z_1& -\bar{z}_2\\
z_2& -\bar{z}_1 \end{pmatrix}\]

For any $N$ we can consider the subgroup $SU(N-1)\subseteq SU(N)$ of elements of the form $1\oplus x$. Denoting by $e_1$ the first basis vector in $\C^N$, we can define a mapping $q:SU(N)\to S^{2N-1}$ by $q(v):=ve_1$. A straightforward calculation shows that $q$ factors over the quotient $SU(N)/SU(N-1)$ and induces a diffeomorphism $SU(N)/SU(N-1)\cong S^{2N-1}$. The function $u$ is in a sense a splitting to $q$:

\begin{prop}
If $\iota:S^{2n-1}\to S^{2^n-1}$ is defined by 
\[\iota(z_1,z_2,\ldots z_n):=
\begin{cases}
(-z_1,z_2,\ldots z_n,0,\ldots,0)\quad \mbox{for $n$ even}\\
(-\bar{z}_1,z_2,\ldots z_n,0,\ldots,0)\quad \mbox{for $n$ odd}\\
\end{cases}\]
and $q:SU(2^{n-1})\to S^{2^n-1}$ is the mapping constructed above, the following identity is satisfied 
\[q\circ u=\iota.\]
\end{prop}

\begin{proof}
We will start with the case when $n$ is even. The first $n$ basis vectors of $S_V^+$ are $1,e_1\wedge e_2, e_1\wedge e_3, \ldots, e_1\wedge e_n$ and 
\[q(u(z))=u(z)1=-z_1+z_2e_1\wedge e_2+z_3e_1\wedge e_3+ \cdots + z_ne_1\wedge e_n.\]
If $n$ is odd, the first basis vectors of $S_V^-$ are $e_1, e_2, \ldots, e_n$. Therefore we have the equality
\[q(u(z))=u(z)e_1=-\bar{z}_1e_1+z_2e_2+\cdots z_ne_n.\]
\end{proof}

Consider $\alpha_+:=\phi_+(\rd z\oplus 0)$ and $\alpha_-:=\phi_+(0\oplus \rd \bar{z})$ as elements in $T^*S^{2n-1}\otimes End_\C(E_n)$.  For an element $\mathbbm{k}=(k_1,\ldots, k_{2l-1})\in \{+,-\}^{2l-1}$ we define $\alpha_\mathbbm{k}:=\alpha_{k_1}\alpha_{k_2}\cdots \alpha_{k_{2l-1}}\in \wedge^{2l+1}T^*S^{2n-1}\otimes End_\C(E_n)$. Define the set $\Gamma_l^+$ as the set of $\mathbbm{k}\in \{+,-\}^{2l-1}$ such that the number of $+$ in $\mathbbm{k}$ is $l$. Similarly $\Gamma_l^-$ is defined as the set of $\mathbbm{k}\in \{+,-\}^{2l-1}$ such that the number of $-$ in $\mathbbm{k}$ is $l$. The number of elements in $\Gamma_l^\pm$ can be calculated as
\[|\Gamma_l^+|=|\Gamma_l^-|=
\begin{pmatrix} 
2l-1\\ 
l-1
\end{pmatrix}=
\frac{(2l-1)!}{l!(l-1)!}.\]

\begin{lem}
\label{evencre}
For any $\mathbbm{k}\in \{+,-\}^{2l-1}$ we have the equalities
\[\tra(z_+\alpha_\mathbbm{k})=
\begin{cases}
0\quad\mbox{if}\quad \mathbbm{k}\notin \Gamma_l^-\\
\;\\
(-1)^n2^{n-1}l! \sum_{m_1,m_2,\ldots m_l} z_{m_1}\rd \bar{z}_{m_1}\bigwedge_{j=2}^l \rd z_{m_j}\wedge\rd \bar{z}_{m_j}\quad\mbox{if}\quad \mathbbm{k}\in \Gamma_l^-\\
\end{cases}\]
\[\tra(\bar{z}_-\alpha_\mathbbm{k})=
\begin{cases}
0\quad\mbox{if}\quad \mathbbm{k}\notin \Gamma_n^+\\
\;\\
(-1)^{n+1}2^{n-1}l!\sum_{m_1,m_2\ldots, m_l} \bar{z}_{m_1}\rd z_{m_1}\bigwedge_{j=2}^l\rd z_{m_j}\wedge\rd \bar{z}_{m_j}\quad\mbox{if}\quad \mathbbm{k}\in \Gamma_l^+\\
\end{cases}\]
Here $\tra$ denotes the matrix trace in $End_\C(E_n)$.
\end{lem}

The proof is a straightforward, but rather lengthy, calculation using the relations in the Clifford algebra, so we omit the proof. We will use the notation $\rd V$ for the normalized volume measure on $S^{2n-1}$:
\begin{align}
\label{snvolume}
\rd V&=\frac{(n-1)!}{2\pi^n}\sum_{k=1}^{2n}(-1)^{k-1}x_k \rd x_1\wedge\cdots \wedge \rd x_{k-1}\wedge \rd x_{k+1}\wedge \cdots\wedge \rd x_{2n}=\\
\label{sncomplexvolume}
&=\frac{(n-1)!}{2(2\pi i)^n}\sum_{k=1}^{n}\bar{z}_k \rd z_k\wedge_{j\neq k}(\rd z_j\wedge\rd \bar{z}_j)-z_k \rd \bar{z}_k\wedge_{j\neq k}(\rd z_j\wedge\rd \bar{z}_j).
\end{align}
That $\rd V$ is normalized follows from that the $2n-1$-form $\omega$ on $S^{2n-1}$, defined by
\[\omega=\sum_{k=1}^{2n}(-1)^{k-1}x_k \rd x_1\wedge\cdots \wedge \rd x_{k-1}\wedge \rd x_{k+1}\wedge \cdots\wedge \rd x_{2n},\]
satisfies that, if we change to spherical coordinates, the form $r^{2n-1}\rd r \wedge \omega$ coincide with the volume form on $\C^n$. Since $\int_\C\e^{-|z|^2}\rd m=\pi$, where $m$ denotes Lebesgue measure, Fubini's Theorem implies that $\int_{\C^n}\e^{-|z|^2}\rd m=\pi^n$ and 
\[\pi^n=\int_{\C^n}\e^{-|z|^2}\rd m=\int_0^\infty \e^{-r^2}r^{2n-1}\rd r\int_{S^{2n-1}} \omega=\frac{(n-1)!}{2} \int_{S^{2n-1}} \omega.\]

Recall that if $g:Y\to \mathrm{GL}_N(\C)$ is a smooth mapping, the Chern-Simons character of $g$ is an element of the odd de Rham cohomology $H^{odd}_{dR}(Y)$ defined as 
\[\ch[g]=\sum _{k=0}^\infty \frac{(k-1)!}{(2\pi i)^k(2k-1)!}\tra (g^{-1}\rd g)^{2k-1}.\]
See more in Chapter $1.8$ in \cite{zhang}. We will denote the $2k-1$-degree term by $\ch_{2k-1}[g]$. The cohomology class of $\ch[g]$ only depends on the homotopy class of $g$ so the Chern-Simons character induces a group homomorphism $\ch:K_1(C^\infty(Y))\to H^{odd}_{dR}(Y)$.

\begin{lem}
\label{genktlem}
The mapping $u$ defined in \eqref{udef} satisfies 
\[\ch[u]=(-1)^n\rd V.\]
\end{lem}

\begin{proof}
Since the odd de Rham cohomology of $S^{2n-1}$ is spanned by the volume form it will be sufficient to show that $\ch_{2n-1}[u]=(-1)^n\rd V$. First we observe the identity $u^*\rd u=-\rd u^* u$, which follows from Proposition \ref{unitlem}. This fact implies 
\[(u^*\rd u)^{2n-1}=(-1)^{n-1}u^*\underbrace{\rd u\,\rd u^*\,\cdots\,\rd u^*\,\rd u}_{2n-1 \quad\mbox{factors}}.\]
Our second observation is 
\[u^*\rd u=-\frac{1}{2}(z+\bar{z})(\rd z+\rd \bar{z}) \quad\mbox{and}\quad \rd u^*\,\rd u=-\frac{1}{2}(\rd z+\rd\bar{z})(\rd z+\rd \bar{z}).\]
Therefore
\[(u^*\rd u)^{2n-1}=-\frac{1}{2^n}(z+\bar{z})(\rd z+\rd \bar{z})^{2n-1}.\]
Because of Lemma \ref{evencre} we have the equalities
\begin{align*}
\tra((z+\bar{z})&(\rd z+\rd \bar{z})^{2n-1})=\sum_{\mathbbm{k}\in \Gamma_n^+}\tra (\bar{z}\alpha_{\mathbbm{k}})+\sum_{\mathbbm{k}\in \Gamma_n^-}\tra (z\alpha_{\mathbbm{k}})=\\
&=\sum_{\mathbbm{k}\in \Gamma_n^+}(-1)^{n+1} 2^{n-1}(n-1)!n! \sum_{k=1}^n  \bar{z}_k\rd z_k\wedge_{j\neq k}(\rd z_j\wedge\rd \bar{z}_j)+\\
&+\sum_{\mathbbm{k}\in \Gamma_n^-}(-1)^n2^{n-1}(n-1)!n! \sum_{k=1}^n z_k\rd \bar{z}_k\wedge_{j\neq k}(\rd z_j\wedge\rd \bar{z}_j)=\\
=(-1)^{n+1}2^{n-1}(2n-& 1)! \sum_{k=1}^n \left(\bar{z}_k\rd z_k\wedge_{j\neq k}(\rd z_j\wedge\rd \bar{z}_j)-z_k\rd \bar{z}_k\wedge_{j\neq k}(\rd z_j\wedge\rd \bar{z}_j)\right)=\\
&=\frac{(-1)^{n+1}2^{n}(2\pi i)^n(2n-1)!}{(n-1)!}\rd V.
\end{align*}
Finally, adding all results together we come to the conclusion of the Lemma: 
\[\tra(u^*\rd u)^{2n-1}=-\frac{1}{2^n}\tra((z+\bar{z})(\rd z+\rd \bar{z})^{2n-1})=(-1)^n \frac{(2\pi i)^n(2n-1)!}{(n-1)!}\rd V.\]

\end{proof}

To generalize the construction of $u$ to an arbitrary manifold we need to cut down $u$ at "infinity". We define the smooth function $\xi_0:[0,\infty)\to \R$ as
\[\xi_0(x):=
\begin{cases} 
\e^{-\frac{4}{x^2}}, \quad x>0\\
0, \;\,\,\quad\quad x=0
\end{cases}\]
and the smooth function $\xi:S^{2n-1}\to \C^n$ by 
\[\xi(z):=\xi_0(|1-\Re(z_1)|)z+\left(\xi_0(|1-\Re(z_1)|)-1,0, 0,\ldots, 0\right).\]
By standard methods it can be proved that for any natural number $k$ and any vector fields $X_1,X_2, \ldots X_l$ on $S^{2n-1}$ the function $\xi$ satisfies 
\begin{align}
\label{conve}
|\xi(z)-(-1,0,\ldots,0)|&=\mathcal{O}(|1-\Re(z_1)|^k)\quad\mbox{and}\\
\label{convetva}
|X_1X_2\cdots X_l\xi(z)|&=\mathcal{O}(|1-\Re(z_1)|^k)\quad\mbox{as}\quad z\to (1,0,\ldots,0).
\end{align}
Furthermore, the length of $\xi(z)$ is given by
\[|\xi(z)|^2=2(\Re(z_1)+1)(\xi_0(|1-\Re(z_1)|)^2-\xi_0(|1-\Re(z_1)|))+1\]
so $|\xi(z)|>0$ for all $z\in S^{2n-1}$.

Using the function $\xi$ we define the smooth function $\tilde{u}:S^{2n-1}\to GL_{2^{n-1}}(\C)$ by 
\[\tilde{u}(z):=\frac{1}{2}(e_{1,+}+e_{1,-})(\xi(z)_++\overline{\xi(z)}_-).\]
The function $\tilde{u}$ is well defined since 
\[\tilde{u}(z)^*\tilde{u}(z)=|\xi(z)|^2>0.\]
Observe that we may express $\tilde{u}$ in terms of $u$ as
\[\tilde{u}(z)=\xi_0(|1-\Re(z_1)|)(u(z)-1)+1.\]

If we choose a diffeomorphism $\tau:B_{2n-1}\cong S^{2n-1}\setminus\{(1,0,\ldots,0)\}$ the equation \eqref{conve} and \eqref{convetva} implies that the function $\tau^*\tilde{u}$ can be considered as a smooth function $B_{2n-1}\to GL_{2^{n-1}}(\C)$ such that $\tau^*\tilde{u}-1$ vanishes to infinite order at the boundary of $B_{2n-1}$. The particular choice of $\tau$ as the stereographic projection
\[\tau(y):=\left(2|y|^2-1,2\sqrt{1-|y|^2}y\right)\]
will give a function $\tau^*\tilde{u}$ of the form
\begin{align*}
\tau^*\tilde{u}(y)&=\e^{-\frac{1}{(1-|y|^2)^2}}(u(\tau(y))-1)+1=\\
&=\frac{\e^{-\frac{1}{(1-|y|^2)^2}}}{2}(e_{1,+}+e_{1,-})(\tau(y)_++\overline{\tau(y)}_-)+1-\e^{-\frac{1}{(1-|y|^2)^2}}.
\end{align*}

\begin{lem}
\label{spherehomotopy}
There is a homotopy of smooth functions $S^{2n-1}\to GL_{2^{n-1}}(\C)$ between $\tilde{u}$ and $u$. Therefore $\ch[\tilde{u}]-\ch[u]$ is an exact form.
\end{lem}

\begin{proof}
We can take the homotopy $w:S^{2n-1}\times [0,1]\to GL_{2^{n-1}}(\C)$ as 
\[w(z,t)=\xi_t(|1-\Re(z_1)|)(u(z)-1)+1,\]
where 
\[\xi_t(x):=\e^{-\frac{4(1-t)}{x^2}}.\]
Clearly, $w:S^{2n-1}\times [0,1]\to GL_{2^{n-1}}(\C)$ is a smooth function and $w(z,0)=\tilde{u}(z)$ and $w(z,1)=u(z)$. 
\end{proof}

In the general case, let $Y$ be a compact, connected, orientable manifold of odd dimension $2n-1$. If we take an open subset $U$ of $Y$ with coordinates $(x_{i})_{i=1}^{2n-1}$ such that 
\[U=\{x:\sum_{i=1}^{2n-1}|x_{i}(x)|^2<1\},\]
the coordinates define a diffeomorphism $\nu:U\cong B_{2n-1}$. We can define the functions $g,\tilde{g}:Y\to GL_{2^{n-1}}(\C)$ by 
\begin{equation}
\label{lipschitzsymbol}
g(x):=\begin{cases}
u(\tau\nu(x))\quad\mbox{for} \quad x\in U\\
1\quad\mbox{for} \quad x\notin U
\end{cases}
\end{equation}
\begin{equation}
\label{localsymbol}
\tilde{g}(x):=\begin{cases}
\tilde{u}(\tau\nu(x))\quad\mbox{for} \quad x\in U\\
1\quad\mbox{for} \quad x\notin U
\end{cases}
\end{equation}
If we let $\tilde{\nu}:Y\to S^{2n-1}$ be the Lipschitz continuous function defined by 
\begin{equation}
\label{nuext}
\tilde{\nu}(x)=
\begin{cases}
\tau(\nu(x)) \quad\mbox{for} \quad x\in U\\
(1,0,\ldots, 0)\quad\mbox{for} \quad x\notin U
\end{cases}
\end{equation}
the functions $\tilde{g}$ and $g$ can be expressed as $g=\tilde{\nu}^*u$ and $\tilde{g}=\tilde{\nu}^*\tilde{u}$. The function $\tilde{g}$ is smooth and the function $g$ is Lipschitz continuous.

\begin{sats}
\label{degsats}
Denoting the normalized volume form on $Y$ by $\rd V_Y$, the function $\tilde{g}$ satisfies  
\begin{equation}
\label{matrixpairing}
\ch[\tilde{g}]=(-1)^n \rd V_Y,
\end{equation}
in $H^{odd}_{dR}(Y)$. Thus, if $f:X\to Y$ is a smooth mapping 
\[\deg(f)=(-1)^n\int_X f^*\ch[\tilde{g}]\]
\end{sats}

\begin{proof}
By Lemma \ref{spherehomotopy} and Lemma \ref{genktlem} we have the identities
\begin{align*}
\int _Y\ch[\tilde{g}]&= \int _{U}\ch_{2n-1}[\tilde{g}]= \int _{U}\tilde{\nu}^*\ch_{2n-1}[\tilde{u}]=\\
&=\int _{S^{2n-1}}\ch_{2n-1}[\tilde{u}]=\int_{S^{2n-1}}\ch_{2n-1}[u]=(-1)^n.
\end{align*}
Therefore we have the identity $\ch_{2n-1}[\tilde{g}]=(-1)^n \rd V_Y$. Since $\ch[\tilde{g}]-\ch_{2n-1}[\tilde{g}]$ is an exact form on $U$ and vanishes to infinite order at $\partial U$ the Theorem follows.
\end{proof}

\section{Toeplitz operators and their index theory}

In this section we will give the basics of integral representations of holomorphic functions and the Henkin-Ramirez integral representation, we will more or less pick out the facts of \cite{forn}, \cite{henkleit} and \cite{range} relevant for our purposes. After that we will review the theory of Toeplitz operators on the Hardy space on the boundary of a strictly pseudo-convex domain. We will let $M$ denote a Stein manifold and we will assume that $\Omega\subseteq M$ is a relatively compact, strictly pseudo-convex domain with smooth boundary. 

Consider the Hilbert space $L^2(\partial\Omega)$, in some Riemannian metric on $\partial\Omega$. We will use the notation $H^2(\partial\Omega)$ for the Hardy space, that is defined as the space of functions in $L^2(\partial\Omega)$ with holomorphic extensions to $\Omega$. The subspace $H^2(\partial\Omega)\subseteq L^2(\partial\Omega)$ is a closed subspace so there exists a unique orthogonal projection $P_{\partial\Omega}:L^2(\partial\Omega)\to H^2(\partial\Omega)$ called the Szeg\"o projection. We will consider the Henkin-Ramirez projection, see \cite{henka}, \cite{ramirein} and the generalization in \cite{henkleit} to Stein manifolds, which we will denote by $P_{HR}:L^2(\partial\Omega)\to H^2(\partial\Omega)$ and call the HR-projection. The HR-projection is not necessarily orthogonal but is often possible to calculate explicitly, see \cite{range},  and easier to estimate. We will briefly review its construction in the case $M=\C^n$ following Chapter VII of \cite{range}. The construction of the HR-projection on a general Stein manifold is somewhat more complicated, but the same estimates hold so we refer the reader to the construction in \cite{henkleit}. 

The kernel of the HR-projection should be thought of as the first terms in a Taylor expansion of the Szeg\"o kernel. This idea is explained in \cite{kerzmanstein}. The HR-kernel contains the most singular part of the Szeg\"o kernel and the HR-kernel can be very explicitly estimated at its singularities. This is our reason to use the HR-projection instead of the Szeg\"o projection. If $\Omega$ is defined by the strictly pluri-subharmonic function $\rho$ a function $\Phi=\Phi(w,z)$ is defined as the smooth global extension of the Levi polynomial 
\[F(w,z):=\sum_{j=1}^n\frac{\partial\rho}{\partial w_j}(w)(w_j-z_j)-\frac{1}{2}\sum_{j,k=1}^n\frac{\partial^2\rho}{\partial w_j\partial w_k}(w)(w_j-z_j)(w_k-z_k)\]
from the diagonal in $\Omega\times \Omega$ to the whole of $\overline{\Omega}\times \overline{\Omega}$, see more in Chapter V$.1.1$ and Chapter VII$.5.1$ of \cite{range}. If we take $c>0$ such that $\partial\bar{\partial}\rho\geq c$ there is an $\epsilon>0$ such that
\begin{equation}
\label{phihatestimate}
2\Re \Phi (w,z)\geq \rho(w)-\rho(z)+c |z-w|^2, \quad \mbox{for}\quad |z-w|<\epsilon,
\end{equation}
see more in equation $1.6$, Chapter V$.1.1$ of \cite{range}. By Lemma $1.5$ of Chapter VII of \cite{range} the function $\Phi$ satisfies the following estimate
\begin{equation}
\label{phiintegral}
\int  _{\partial\Omega}\frac{\rd V(w)}{|\Phi (w,z)|^{n+\beta}}\lesssim 1
\end{equation}
where $\rd V$ denotes the volume measure on $\partial\Omega$ if $\beta<0$ and a similar estimate with the roles of $z$ and $w$ interchanged.  Here we used the standard notation $a\lesssim b$ for the statement that there exists a constant $C>0$ such that $a\leq Cb$.

By Theorem $3.6$, Chapter VII of \cite{range} we can associate with $\Phi$ a function $H_{\partial\Omega}$ in $\Omega\times \Omega$ holomorphic in its second variable such that if $g\in L^1(\Omega)$ is holomorphic it has the integral representation: 
\[g(z)=\int_{\partial\Omega} H_{\partial\Omega}(w,z)f(w)\rd V(w).\]
For the function $H_{\partial\Omega}$ the estimate
\begin{equation}
\label{ginphi}
|H_{\partial\Omega}(z,w)|\lesssim |\Phi(w,z)| ^{-n},
\end{equation}
holds in $\partial\Omega\times \partial\Omega$, see more in Proposition $3.1$, Chapter VII of \cite{range}. Since $\Phi$ satisfies the estimate \eqref{phihatestimate} where $c$ is the infimum of $\partial \bar{\partial}\rho$ the construction of a HR-projection does give an $L^2$-bounded projection for strictly pseudo-convex domains. If $\Omega$ is weakly pseudo-convex the situation is more problematic and not that well understood partly due to problems estimating solutions to the $\bar{\partial}$-equation in weakly pseudo-convex domains. By Proposition $3.8$ of Chapter VII$.3.1$ in \cite{range} the kernel $H_{\partial\Omega}$ satisfies the estimate 
\begin{equation}
\label{diffhenkhenk}
|H_{\partial\Omega}(z,w)-\overline{H_{\partial\Omega}(w,z)}|\lesssim |\Phi (z,w)|^{-n+1/2}.
\end{equation}
The estimate \eqref{diffhenkhenk} will be crucial when proving that $P_{\partial\Omega}-P_{HR}$ is in the Schatten class. The kernel $H_{\partial\Omega}$ determines a bounded operator $P_{HR}$ on $L^2(\partial\Omega)$ by Theorem $3.6$ of Chapter $VII.3$ in \cite{range}. Since the range of $P_{HR}$ is contained in $H^2(\partial\Omega)$ and $g=P_{HR}g$ for any $g\in H^2(\partial\Omega)$ it follows that $P_{HR}:L^2(\partial\Omega)\to H^2(\partial\Omega)$ is a projection.

We will now present some facts about Toeplitz operators on the Hardy space of a relatively compact strictly pseudo-convex domain $\Omega$ in a complex manifold $M$. Our operators are associated with the Szeg\"o projection since the theory becomes somewhat more complicated when a non-orthogonal projection is involved. For any dimension $N$ we denote by $C(\partial\Omega,M_N)$ the $C^*$-algebra of continuous functions $\partial\Omega \to M_N$, the algebra of complex $N\times N$-matrices. The algebra $C(\partial\Omega ,M_N)$ has a representation $\pi:C(\partial\Omega ,M_N)\to \Bo(L^2(\partial\Omega)\otimes \C^N)$ which is given by pointwise multiplication. We define the linear mapping 
\[T:C(\partial\Omega ,M_N)\to \Bo(H^2(\partial\Omega)\otimes \C^N), \quad a\mapsto P_{\partial\Omega} \pi(a) P_{\partial\Omega}.\] 
Here we identify $P_{\partial\Omega}$ with the projection $L^2(\partial\Omega)\otimes \C^N\to H^2(\partial\Omega)\otimes \C^N$. An operator of the form $T(a)$ is called a Toeplitz operator on $\partial\Omega$. Toeplitz operators are well studied, see for instance \cite{boutetdemonvel}, \cite{connes}, \cite{guehig} and \cite{grigone}. The representation $\pi$ satisfies $[P_{\partial\Omega},\pi(a)]\in \Ko(L^2(\partial\Omega)\otimes \C^N)$ for any $a\in C(\partial\Omega ,M_N)$, see for instance \cite{boutetdemonvel} or Theorem \ref{sumsats} below. Here we use the symbol $\Ko$ to denote the algebra of compact operators. The fact that $P_{\partial\Omega}$ commutes with continuous functions up to a compact operator implies the property 
\begin{equation}
\label{almost}
T(ab)-T(a)T(b)\in \Ko(H^2(\partial\Omega)\otimes \C^N).
\end{equation}
Furthermore, $T(a)$ is compact if and only if $a=0$. Let us denote the Calkin algebra $\Bo(\He)/\Ko(\He)$ by $\mathcal{C}(\He)$ and the quotient mapping $\Bo(\He)\to \mathcal{C}(\He)$ by $\mathfrak{q}$. Equation \eqref{almost} implies that the mapping
\[\beta:=\mathfrak{q}\circ T: C(\partial\Omega ,M_N)\to \mathcal{C}(H^2(\partial\Omega)\otimes \C^N)\]
is an injective $*$-homomorphism. 

By the Boutet de Monvel index formula, from \cite{boutetdemonvel}, if the symbol $a$ is invertible and smooth the index of the Toeplitz operator $T(a)$ has the analytic expression:
\begin{equation}
\label{indform}
\ind (T(a))= -\int_{\partial\Omega}\ch[a]\wedge Td(\Omega),
\end{equation}
see more in Theorem $1$ in \cite{boutetdemonvel}, and the remarks thereafter. The mapping $a\mapsto \ind (T(a))$, defined on functions $a:\partial \Omega\to \mathrm{GL}_N(\C)$ is homotopy invariant, so it extends to a mapping $\ind:K_1(C^\infty(\partial\Omega))\to \Z$. Here $K_1(C^\infty(\partial\Omega))$ denotes the odd K-theory of the Frechet algebra $C^\infty(\partial\Omega)$ which is defined as homotopy classes of invertible matrices with coefficients in $C^\infty(\partial\Omega)$, see more in \cite{blacker}. 

\begin{sats}
\label{indexdegreeform}
Suppose that $\Omega\subseteq M$ is a relatively compact strictly pseudo-convex bounded domain with smooth boundary, $Y$ is a compact, orientable manifold of dimension $2n-1$ and $g:Y\to GL_{2^{n-1}}(\C)$ is the mapping defined in \eqref{localsymbol}. If $f:\partial \Omega\to Y$ is a continuous function, then 
\begin{equation}
\label{indexdegree}
\deg(f)=(-1)^{n+1}\ind (P_{\partial\Omega} \pi(g\circ f)P_{\partial\Omega}).
\end{equation}
\end{sats}

\begin{proof}
If we assume that $f$ is smooth, the index formula of Boutet de Monvel, see above in equation \eqref{indform}, implies that the index of $P_{\partial\Omega} \pi(g\circ f)P_{\partial\Omega}$ satisfies 
\[\ind (P_{\partial\Omega} \pi(g\circ f)P_{\partial\Omega})=-\int _{\partial \Omega}f^*\ch[\tilde{g}] \wedge Td(\Omega)=-\int _{\partial \Omega}f^*\ch[\tilde{g}] =(-1)^{n+1}\deg(f),\]
where the first equality follows from $g$ and $\tilde{g}$ being homotopic, see Lemma \ref{spherehomotopy}, and the last two equalities follows from Theorem \ref{degsats}. The general case follows from the fact that both hand sides of \eqref{indexdegree} is homotopy invariant.
\end{proof}

Theorem \ref{indexdegreeform} does in some cases hold with even looser regularity conditions on $f$. Since both sides of the equation \eqref{indexdegree} are homotopy invariants the Theorem holds for any class of functions which are homotopic to smooth functions in such sense that both sides in \eqref{indexdegree} are well defined and depend continuously on the function. For instance, if $\Omega$ is a bounded symmetric domain we may take $f:\partial\Omega\to Y$ to be in the $VMO$-class. It follows from \cite{bergcobuzhu} that if $w:\partial\Omega\to \mathrm{GL}_N$ has vanishing mean oscillation and $\Omega$ is a bounded symmetric domain, the operator $P_{\partial\Omega} wP_{\partial\Omega}$ is Fredholm. By \cite{brenir} the degree of a $VMO$-function is well defined and depends continuously on $f$ without any restriction on the geometry. To be more precise, there is a one-parameter family $(f_t)_{t\in (0,1)}\subseteq C(\partial\Omega,Y)$ such that $f_t\to f$ in $VMO$ when $t\to 0$ and $\deg(f)$ is defined as $\deg(f_t)$ for $t$ small enough. Since the index of a Fredholm operator is homotopy invariant the degree of a function $f:\partial\Omega\to Y$ in $VMO$ satisfies 
\[\deg f=(-1)^{n+1}\ind (P_{\partial\Omega} \pi(g\circ f_t) P_{\partial\Omega})=(-1)^{n+1}\ind (P_{\partial\Omega} \pi(g\circ f) P_{\partial\Omega}).\]
\[\]

Our next task will be calculating the index of Toeplitz operators with non-smooth symbol. For $p\geq 1$, let $\ellL^p(\He)\subseteq \Bo(\He)$ denote the ideal of Schatten class operators on a separable Hilbert space $\He$, so $T\in \ellL^p(\He)$ if and only if $\tra((T^*T)^{p/2})<\infty$. An exact description of integral operators belonging to this class exists only for $p=2$. However, for $p>2$ there exists a convenient sufficient condition on the kernel, found in \cite{russo}. We will return to this subject a little later. Suppose that $\pi:\mathcal{A}\to \Bo(\He)$ is a representation of a $\C$-algebra $\mathcal{A}$ and $P$ is a projection such that $[P,\pi(a)]\in \ellL^{p}(\He)$ for all $a\in \mathcal{A}$ and $P-P^*\in \ellL^{p}(\He)$. Atkinson's Theorem implies that if $a$ is invertible, $P\pi(a)P$ is Fredholm. The operator $F:=2P-1$ has the properties 
\begin{equation}
\label{fredprop}
F^2=1\quad\mbox{and} \quad F-F^*,\;[F,\pi(a)]\in \ellL^{p}(\He).
\end{equation}
If $\pi$ and $F$ satisfy the conditions in equation \eqref{fredprop} the pair $(\pi,F)$ is called a $p$-summable odd Fredholm module. If the pair $(\pi,F)$ satisfies the requirement in equation \eqref{fredprop} but with $\ellL^p(\He)$ replaced by $\Ko(\He)$ the pair $(\pi,F)$ is a bounded odd Fredholm module. For a more thorough presentation of Fredholm modules, e.g. Chapter VII and VIII of \cite{blacker}. Since our focus is on Toeplitz operators we will call $(\pi,P)$ a \emph{Toeplitz pair} if $(\pi,2P-1)$ is a bounded odd Fredholm module and $(\pi,P)$ is said to be $p$-summable if $(\pi,2P-1)$ is.

The condition that $L:=P^*-P\in \ellL^{p}(\He)$ can be interpreted in terms of the orthogonal projection $\tilde{P}$ to the Hilbert space $P\He$. Using that $\tilde{P}P=P$ and $P\tilde{P}=\tilde{P}$ we obtain the identity
\begin{equation}
\label{projdiff}
\tilde{P}L=\tilde{P}P^*-\tilde{P}P=\tilde{P}-P.
\end{equation}
Thus the condition $P^*-P\in \ellL^p(\He)$ is equivalent to the property $\tilde{P}-P\in \ellL^p(\He)$.

A Toeplitz pair $(\pi,P)$ over a topological algebra $\mathcal{A}$ defines a mapping $a\mapsto \ind(P\pi(a)P)$ on the invertible elements of $\mathcal{A}\otimes M_N$ for any $N$. Since the index is homotopy invariant, the association $a\mapsto \ind(P\pi(a)P)$ induces the mapping $\ind:K_1(\mathcal{A})\to \Z$, where $K_1(\mathcal{A})$ denotes the odd $K$-theory of $\mathcal{A}$, see \cite{blacker}. 

A. Connes placed the index theory for $p$-summable Toeplitz pairs in a suitable homological picture using cyclic homology in \cite{connesncdg}. We will consider Connes' original definition of cyclic cohomology which simplifies the construction of the Chern-Connes character. The notation $\mathcal{A}^{\otimes k}$ will be used for the $k$-th tensor power of $\mathcal{A}$. The Hochschild differential $b:\mathcal{A}^{\otimes k}\to \mathcal{A}^{\otimes k-1}$ is defined as 
\begin{align*}
b(x_0\otimes x_1\otimes \cdots\otimes x_k\otimes x_{k+1})&:=(-1)^{k+1}x_{k+1}x_0\otimes x_1\otimes \cdots \otimes x_k+\\
+\sum_{j=0}^{k}(-1)^j&x_0\otimes  \cdots \otimes  x_{j-1}\otimes  x_jx_{j+1}\otimes x_{j+2}\otimes \cdots \otimes x_{k+1}.
\end{align*}
The cyclic permutation operator $\lambda:\mathcal{A}^{\otimes k}\to \mathcal{A}^{\otimes k}$ is defined by
\[\lambda(x_0\otimes x_1\otimes \cdot \otimes x_k)=(-1)^kx_k\otimes x_0\otimes \cdots \otimes x_{k-1}.\]
The complex $C^k_\lambda(\mathcal{A})$ is defined as the space of continuous linear functionals $\mu$ on $\mathcal{A}^{\otimes k+1}$ such that $\mu\circ\lambda=\mu$. The Hochschild coboundary operator $\mu\mapsto \mu\circ b$ makes $C^*_\lambda(\mathcal{A})$ into a complex. The cohomology of the complex $C^*_\lambda(\mathcal{A})$ will be denoted by $HC^*(\mathcal{A})$ and is called the cyclic cohomology of $\mathcal{A}$. There is a filtration on cyclic cohomology coming from a linear mapping $S:HC^k(\mathcal{A})\to HC^{k+2}(\mathcal{A})$ which is called the suspension operator or the periodicity operator. For a definition of the periodicity operator, see \cite{connes}. 

The additive pairing between $HC^{2k+1}(\mathcal{A})$ and the odd $K$-theory $K_1(\mathcal{A})$ is defined by 
\[\langle \mu, u\rangle_k:=d_k\;(\mu\otimes \tra)\left(\underbrace{(u^{-1}-1)\otimes(u-1)\otimes \cdots \otimes (u^{-1}-1)\otimes (u-1)}_{2k+2 \quad\mbox{factors}}\right)\]
where we choose the same normalization constant $d_k$ as in Proposition $3$ of Chapter III$.3$ of \cite{connes}:
\[d_k:=\frac{2^{-(2k+1)}}{\sqrt{2i}}\Gamma\left(\frac{2k+3}{2}\right)^{-1}.\]
The choice of normalization implies that for a cohomology class in $HC^{2k+1}(\mathcal{A})$ represented by the cyclic cocycle $\mu$, the pairing satisfies
\[\langle S\mu, u\rangle_{k+1}=\langle \mu, u\rangle_{k},\]
see Proposition $3$ in Chapter III$.3$ of \cite{connes}. Following Definition $3$ of Chapter IV$.1$ of \cite{connes} we define the Connes-Chern character of a $p$-summable Toeplitz pair as the cyclic cocycle: 
\[\che_{2k+1}(\pi,P)(a_0,a_1,\ldots, a_{2k+1}):=c_k\tra (\pi(a_0)[P,\pi(a_1)]\cdots [P,\pi(a_{2k+1})]),\]
for $2k+1\geq p$ where 
\[c_k:=-\sqrt{2i}2^{2k+1}\Gamma\left(\frac{2k+3}{2}\right).\]
This choice of normalization constant implies that 
\[S\che_{2k+1}(\pi,P)=\che_{2k+3}(\pi,P),\]
by Proposition $2$ in Chapter IV$.1$ of \cite{connes}.

\begin{sats}[Proposition $4$ of Chapter IV$.1$ of \cite{connes}]
\label{connesindex}
If $(\pi,P)$ is a $p$-summable Toeplitz pair over $\mathcal{A}$, $2k+1\geq p$ and $a$ is invertible in $\mathcal{A}\otimes M_N$ the index of $P\pi(a)P:P\He\otimes \C^N\to P\He\otimes \C^N$ may be expressed as
\begin{align*}
\ind(P\pi(a)P)&=\langle \che_{2k+1}(\pi,P),a\rangle_k=\\
&=-\tra \left(\pi(a^{-1})[P,\pi(a)][P,\pi(a^{-1})]\cdots [P,\pi(a^{-1})][P,\pi(a)]\right)=\\
&=-\tra(P-\pi(a^{-1})P\pi(a))^{2k+1}.
\end{align*}
\end{sats}

The role of the periodicity operator $S$ in the context of index theory is to extend index formulas to larger algebras. Suppose that $\mu$ is a cyclic $k$-cocycle on an algebra $\mathcal{A}$ which is a dense $*$-subalgebra of a $C^*$-algebra $A$. As is explained in \cite{connes} for functions on $S^1$ and in \cite{grigtwo} for operator valued symbols, the cyclic $k+2m$-cocycle $S^m\mu$ can be extended to a cyclic cocycle on a larger $*$-subalgebra $\mathcal{A}\subseteq \mathcal{A}'\subseteq A$. When $\mu$ is the cyclic cocycle $f_0\otimes f_1\mapsto \int f_0\rd f_1$ on $C^\infty(S^1)$, the $2m+1$-cocycle $S^m\mu$ extends to $C^\alpha(S^1)$ whenever $\alpha(2m+1)>1$ by Proposition $3$ in Chapter III$2.\alpha$ of \cite{connes} and a formula for $S^m\mu$ is given above in \eqref{connesalpha}. Cyclic cocycles of the form $\mu=\che(\pi,P)$ appear in index theory and the periodicity operator can be used to extend index formulas to larger algebras.

The index formula of Theorem \ref{connesindex} holds for Toeplitz operators under a Schatten class condition and to deal with this condition we will need the following theorem of Russo \cite{russo} to give a sufficient condition on an integral operator for it to be Schatten class. Let $X$ denote a $\sigma$-finite measure space. As in \cite{bepa}, for numbers $1\leq p,q <\infty$, the mixed $(p,q)$-norm of a function $k:X\times X\to \C$ is defined by 
\[\|k\|_{p,q}:=\left( \int _X\left( \int_X|k(x,y)|^p\rd x \right) ^{\frac{q}{p}}\rd y\right) ^{\frac{1}{q}}.\]
We denote the space of measurable functions $k:X\times X\to \C$ with finite mixed $(p,q)$-norm by $L^{(p,q)}(X\times X)$. By Theorem $4.1$ of \cite{bepa} the space $L^{(p,q)}(X\times X)$ becomes a Banach space in the mixed $(p,q)$-norm which is reflexive if $1<p,q<\infty$. 

The hermitian conjugate of the function $k$ is defined by $k^*(x,y):=\overline{k(y,x)}$. Clearly, if a bounded operator $K$ has integral kernel $k$, the hermitian conjugate $K^*$ has integral kernel $k^*$. 

\begin{sats}[Theorem $1$ in \cite{russo}]
\label{russothm}
Suppose that $K:L^2(X)\to L^2(X)$ is a bounded operator given by an integral kernel $k$. If $2<p<\infty$ 
\begin{equation}
\label{russoest}
\|K\|_{\ellL^p(L^2(X))}\leq (\|k\|_{p',p}\|k^*\|_{p',p})^{1/2},
\end{equation}
where $p'=p/(p-1)$.
\end{sats}

In the statement of the Theorem in \cite{russo}, the assumption $k\in L^2(X\times X)$ is made. This assumption implies that $K$ is Hilbert-Schmidt and $K\in \ellL^p(L^2(X))$ for all $p>2$ so for our purposes it is not interesting. But since $L^2$-kernels are dense in $L^{(p,q)}$, the non-commutative Fatou lemma, see Theorem $2.7d$ of \cite{simon}, implies \eqref{russoest} for any $k$ for which the right hand side of \eqref{russoest} is finite. Using Theorem \ref{russothm}, we obtain the following formula for the trace of the product of integral operators:

\begin{sats}
\label{russotracethm}
Suppose that $K_j:L^2(X)\to L^2(X)$ are operators with integral kernels $k_j$ for $j=1,\ldots, m$ such that $\|k_j\|_{p',p},\|k_j^*\|_{p',p}<\infty$ for certain $p>2$. Then for $m\geq p$ the operator $K_1K_2\cdots K_m$ is a trace class operator and we have the trace formula
\[\tra(K_1K_2\cdots K_m)=\int _{X^m}\left(\prod_{j=1}^{m}k_{j}(x_j,x_{j+1})\right)\rd x_1\rd x_2\cdots \rd x_m,\]
where we identify $x_{m+1}$ with $x_1$.
\end{sats}

\begin{proof}
The case $p=m=2$ follows if for any $k_1,k_2\in L^2(X\times X)$ we have the trace formula 
\[\tra(KL^*)=\int _{X\times X} k(x,y)\overline{l(x,y)}\rd x\rd y.\]
Consider the sesquilinear form on $\ellL^2(L^2(X))$ defined by
\[ (K,L):=\tra(KL^*)-\int k(x,y)\overline{l(x,y)}\rd x\rd y.\]
Since $\tra(K^*K)=\int_{X\times X}|k(x,y)|^2\rd x\rd y$ the sesquilinear form satisfies $(K,K)=0$ and the polarization identity implies $(K,L)=0$ for any $K,L\in \ellL^2(L^2(X))$.

If the operators $K_j:L^2(X)\to L^2(X)$ are Hilbert-Schmidt, or equivalently they satisfy $k_j\in L^2(X\times X)$, we may take $K=K_1$ and $L^*=K_2K_3\cdots K_m$ so the case $p=m=2$ implies that the operators $K_1,K_2,\ldots, K_m$  satisfy the statement of the Theorem. In the general case,  the Theorem follows from the non-commutative Fatou lemma, see Theorem $2.7d$ of \cite{simon}, since $\ellL^2$ is dense in $\ellL^p$ for $p>2$. 
\end{proof}

\section{The Toeplitz pair on the Hardy space}

As explained in section $2$, for the representation $\pi:C(\partial\Omega )\to \Bo(L^2(\partial\Omega))$ and the Szeg\"o projection $P_{\partial\Omega}$ the commutator $[P_{\partial\Omega},\pi(a)]$ is compact for any continuous $a$. Thus $(\pi,P_{\partial\Omega})$ is a Toeplitz pair over $C(\partial\Omega )$. To enable the use of the index theory of \cite{connes} we will show that the Toeplitz pair $(\pi,P_{\partial\Omega})$ restricted to the subalgebra of H\"older continuous functions $C^\alpha(\partial\Omega)\subseteq C(\partial\Omega)$ becomes $p$-summable. These results will give us analytic degree formulas for H\"older continuous mappings.

\begin{sats}
\label{sumsats}
If $\Omega$ is a relatively compact strictly pseudo-convex domain in a Stein manifold of complex dimension $n$ and $P$ denotes either $P_{HR}$ or $P_{\partial\Omega}$  the operator $[P,\pi(a)]$ belongs to $\ellL^{p}(L^2(\partial\Omega))$ for $a\in C^\alpha(\partial\Omega)$ and for all $p>2n/\alpha$. 
\end{sats}

The proof will be based on Theorem \ref{russothm}. We will start our proof of Theorem \ref{sumsats} by some elementary estimates. We define the measurable function $k_\alpha:\partial\Omega\times\partial \Omega \to \C$ by 
\[k_{\alpha}(z,w):=\frac{|z-w|^{\alpha}}{|\Phi (w,z)|^{n}}.\]

\begin{lem}
\label{pointwiseaest}
The function $k_{\alpha}$ satisfies 
\[k_{\alpha}(z,w)\lesssim |\Phi (w,z)|^{-(n-\frac{\alpha}{2})}\]
for $|z-w|<\epsilon$.
\end{lem}

\begin{proof}
By \eqref{phihatestimate} we have the estimate
\[|z-w|^{\alpha}\lesssim |\Phi (w,z)|^{\alpha/2}.\]
From this estimate the Lemma follows.
\end{proof}

We will use the notation $\rd V$ for the volume measure on $\partial\Omega$.

\begin{lem}
\label{ballgenest}
The function $k_\alpha$ satisfies 
\[\int_{\partial\Omega}|k_\alpha(z,w)|^{p'}\rd V(z)\lesssim 1\]
\[\int_{\partial\Omega}|k_\alpha(z,w)|^{p'}\rd V(w)\lesssim 1\]
whenever
\[(2n-\alpha)p'<2n.\]
\end{lem}

\begin{proof}
We will only prove the first of the estimates in the Lemma. The proof of the second estimate goes analogously. Using \eqref{phihatestimate} for $\Phi $, we obtain 
\[\int_{\partial\Omega}|k_\alpha(z,w)|^{p'}\rd V(z)\lesssim \int_{B_r(w)}|k_\alpha(z,w)|^{p'}\rd V(z),\]
since the function $\Phi $ satisfies $|\Phi (w,z)|>r^2$ outside $B_r (w)$. By Lemma \ref{pointwiseaest} we can estimate the kernel pointwise by $\Phi $ so \eqref{phiintegral} implies 
\[ \int_{B_r(w)}|k_\alpha(z,w)|^{p'}\rd V(z)\lesssim  \int_{B_r(w)}|\Phi (w,z)|^{-p'(n-\frac{\alpha}{2})}\rd V(z)\lesssim 1\]
if $(n-\frac{\alpha}{2})p'<n$.
\end{proof}

\begin{lem}
\label{fundestmix}
The function $k_\alpha$ satisfies $\|k_\alpha\|_{p',p}<\infty$ and $\|k_\alpha^*\|_{p',p}<\infty$ for $p>2n/\alpha$.
\end{lem}

\begin{proof}
By the first estimate in Lemma \ref{ballgenest} we can estimate the mixed norms of $k_\alpha$ as 
\[\|k_\alpha\|_{p',p}^p\lesssim 1,\]
whenever $(2n-\alpha)p'<2n$. The statement $(2n-\alpha)p'<2n$ is equivalent to 
\[\frac{1}{p}=1-\frac{1}{p'}<\frac{\alpha}{2n}\]
which is equivalent to $p>2n/\alpha$. Similarly, the second estimate in Lemma \ref{ballgenest} implies $\|k_\alpha^*\|_{p',p}<\infty$ under the same condition on $p$.
\end{proof}

\begin{lem}
\label{elemest}
Suppose that $a\in C^\alpha(\partial\Omega)$ and let $\kappa_a$ denote the integral kernel of $[P_{HR},\pi(a)]$. The kernel $\kappa_a$ satisfies 
\begin{equation}
\label{ballest}
|\kappa_a(z,w)|\leq \|a\|_{C^\alpha(\partial\Omega)}|k_\alpha(z,w)|,
\end{equation}
where $ \|\cdot\|_{C^\alpha(\partial\Omega)}$ denotes the usual norm in $C^\alpha(\partial\Omega)$.
\end{lem}

\begin{proof}
The integral kernel of $[P_{HR},\pi(a)]$ is given by 
\[\kappa_a(z,w)=(a(z)-a(w))H_{\partial\Omega}(w,z).\]
Since $a$ is H\"older continuous and $H_{\partial\Omega}$ satisfies equation \eqref{ginphi} the estimate \eqref{ballest} follows. 
\end{proof}

\begin{lem}
\label{henkbergdiff}
The HR-projection $P_{HR}$ satisfies $P_{HR}-P_{HR}^*\in \ellL^q(L^2(\partial\Omega))$ for any $q>2n$. Therefore $P_{HR}-P_{\partial\Omega}\in \ellL^q(L^2(\partial\Omega))$ for any $q>2n$.
\end{lem}

\begin{proof}
Let us denote the kernel of the operator $P_{HR}-P_{HR}^*$ by $b$. By \eqref{diffhenkhenk} we have the pointwise estimate $|b(z,w)|\lesssim |\Phi (w,z)|^{-n+1/2}$. Applying Lemma \ref{fundestmix} with $\alpha=0$ and $p'$ such that $(n-1/2)q'=np'$ we obtain the inequality $\|b\|_{q',q}<\infty$ for any $q>2n$. The fact that $P_{HR}-P_{\partial\Omega}\in \ellL^q(L^2(\partial\Omega))$ follows now from \eqref{projdiff}.
\end{proof}

\begin{proof}[Proof of Theorem \ref{sumsats}]
By Lemma \ref{elemest} the integral kernel $\kappa_a$ of $[P_{HR},\pi(a)]$ satisfies $|\kappa_a|\leq \|a\|_{C^\alpha(\partial\Omega)}k_\alpha$. Theorem \ref{russothm} implies the estimate 
\[\|[P_{HR},\pi(a)]\|_{\ellL^p(L^2(\partial\Omega))}\leq \|a\|_{C^\alpha(\partial\Omega)}(\|k_\alpha\|_{p',p}\|k_\alpha^*\|_{p',p})^{1/2}.\]
By Lemma \ref{fundestmix}, $\|k_\alpha\|_{p',p},\|k_\alpha^*\|_{p',p}<\infty$ for $p>2n/\alpha$ so $[P_{HR},\pi(a)]\in \ellL^p(L^2(\partial\Omega))$ for $p>2n/\alpha$. By Lemma \ref{henkbergdiff}, $P_{HR}-P_{\partial\Omega} \in \ellL^{p}(L^2(\partial\Omega))$, so 
\[[P_\Omega,\pi(a)]=[P_{HR},\pi(a)]+[P_\Omega-P_{HR},\pi(a)]\in \ellL^p(L^2(\partial\Omega))\]
for $p>2n/\alpha$ and the proof of the Theorem is complete.
\end{proof}

\section{The index- and degree formula for boundaries of strictly pseudo-convex domains in Stein manifolds}

We may now combine our results on summability of the Toeplitz pairs $(P_{HR},\pi)$ and $(P_{\partial\Omega},\pi)$ into index theorems and degree formulas. The index formula will be proved using the index formula of Connes, see Theorem \ref{connesindex}. 

\begin{sats}
\label{holderindex}
Suppose that $\Omega$ is a relatively compact strictly pseudo-convex domain with smooth boundary in a Stein manifold of complex dimension $n$  and denote the corresponding HR-kernel by $H_{\partial\Omega}$ and the Szeg\"o kernel by $C_{\partial\Omega}$. If $a:\partial\Omega\to\mathrm{GL}_N$ is H\"older continuous with exponent $\alpha$, then for $2k+1>2n/\alpha$ the index formulas hold
\begin{align}
\ind &(P_{\partial\Omega} \pi(a)P_{\partial\Omega})=\ind(P_{HR}\pi(a)P_{HR})=\\
\label{holdindex}&=-\int_{\partial\Omega^{2k+1}} \tra\left(\prod_{j=0}^{2k}(1-a(z_{j-1})^{-1}a(z_j))H_{\partial\Omega}(z_{j-1},z_j)\right)\rd V=\\
\label{holdindexszeg}
&=-\int_{\partial\Omega^{2k+1}} \tra\left(\prod_{j=0}^{2k}(1-a(z_{j-1})^{-1}a(z_j))C_{\partial\Omega}(z_{j-1},z_j)\right)\rd V,
\end{align}
where the integrals in \eqref{holdindex} and \eqref{holdindexszeg} converge.
\end{sats}

\begin{proof}
By Theorem \ref{connesindex} we have 
\[\ind (P_{\partial\Omega} \pi(a)P_{\partial\Omega})=-\tra(P_{\partial\Omega} -\pi(a^{-1})P_{\partial\Omega}\pi(a))^{2k+1}\]
and by Theorem \ref{russotracethm} the trace has the form 
\begin{align*}
-\tra(P_{\partial\Omega} -\pi(a^{-1})&P_{\partial\Omega}\pi(a))^{2k+1}=\\
&=-\int_{\partial\Omega^{2k+1}} \tra\left(\prod_{j=0}^{2k}(1-a(z_{j-1})^{-1}a(z_j))C_{\partial\Omega}(z_{j-1},z_j)\right)\rd V.
\end{align*}
Similarly, the index for $P_{HR}\pi(a)P_{HR}$ is calculated. The Theorem follows from the identity $\ind (P_{\partial\Omega} \pi(a)P_{\partial\Omega})=\ind(P_{HR}\pi(a)P_{HR})$ since Lemma \ref{henkbergdiff} implies that $P_{\partial\Omega} \pi(a)P_{\partial\Omega}-P_{HR}\pi(a)P_{HR}$ is compact.
\end{proof}

Theorem \ref{holderindex} has an interpretation in terms of cyclic cohomology. Define the cyclic $2n-1$-cocycle $\chi_{\partial\Omega}$ on $C^\infty(\partial\Omega)$ by 
\[\chi_{\partial\Omega}:=\sum_{k=0}^n S^k \omega_k,\]
where $\omega_k$ denotes the cyclic $2n-2k-1$-cocycle given by the Todd class $Td_k(\Omega)$ in degree $2k$ as
\[\omega_k(a_0,a_1,\ldots, a_{2n-2k-1}):=\int _{\partial\Omega}a_0 \rd a_1\wedge \rd a_2\wedge \cdots \wedge \rd a_{2n-2k-1}\wedge Td_k(\Omega).\]
Similarly to Proposition $13$, Chapter III.$3$ of \cite{connes}, we have the following:

\begin{sats}
\label{cycref}
The cyclic cocycle $S^{m}\chi_{\partial\Omega}$ defines the same cyclic cohomology class on $C^\infty(\partial\Omega)$ as 
\begin{align*}
\tilde{\chi}_{\partial\Omega}&(a_0,a_1,\ldots, a_{2n+2m-1}):=\\
&:=\int_{\partial\Omega^{2n+2m-1}} \tra\left(a_0(z_0)\prod_{j=1}^{2n+2m-1}(a_j(z_{j})-a_j(z_{j-1}))C_{\partial\Omega}(z_{j-1},z_j)\right)\rd V,
\end{align*}
where we identify $z_{2n+2m-1}=z_0$. Furthermore, the cyclic cocycle $\tilde{\chi}_{\partial\Omega}$ extends to a cyclic $2n+2m-1$-cocycle on $C^\alpha(\partial\Omega)$ if $m>(2n(1-\alpha)+\alpha)/2\alpha$.
\end{sats}

Returning to the degree calculations, to express the degree of a H\"older continuous function we will use Theorem \ref{indexdegreeform} and Theorem \ref{holderindex}. In order to express the formulas in Theorem \ref{holderindex} directly in terms of $f$ we will need some notations. Let $\langle\cdot,\cdot\rangle$ denote the scalar product on $\C^n$. The symmetric group on $m$ elements will be denoted by $S_m$. We will consider $S_m$ as the group of bijections on the set $\{1,2,\ldots, m\}$ and identify the element $m+1$ with $1$ in the set $\{1,2,\ldots, m\}$. 

For $2l\leq m$ we will define a function $\epsilon_l:S_m\to \{0,1,-1\}$ which we will refer to the order parity. If $\sigma\in S_m$ satisfies that there is an $i\in \{\sigma(1),\sigma(2),\ldots \sigma(2l-1),\sigma(2l)\}$ such that $i+1,i-1\notin \{\sigma(1),\sigma(2),\ldots \sigma(2l-1),\sigma(2l)\}$ we set $\epsilon_l(\sigma)=0$. If $\sigma$ does not satisfy this condition the order parity of $\sigma$ is set as $(-1)^k$, where $k$ is the smallest number of transpositions needed to map the set $\{\sigma(1),\sigma(2),\ldots \sigma(2l-1),\sigma(2l)\}$, with $j$ identified with $j+m$, to a set of the form $\{j_1,j_1+1,j_2,j_2+1,\ldots , j_{l},j_l+1\}$ where $1\leq j_1<j_2<\cdots <j_{l}\leq m$.

\begin{prop}
\label{nsch}
The function $u$ satisfies
\begin{align*}
\tra&\left(\prod_{i=0}^{2k}(1-u(z_{i-1})^*u(z_{i}))\right)=\\
&\quad=\sum _{l=0}^{2k+1}\sum_{\sigma\in S_{2(2k+1)}} (-1)^l2^{n-l-1}\epsilon_l(\sigma)\langle z_{\sigma(1)},z_{\sigma(2)}\rangle\langle z_{\sigma(3)},z_{\sigma(4)}\rangle\cdots \langle z_{\sigma(2l-1)},z_{\sigma(2l)}\rangle,
\end{align*}
where we identify $z_m$ with $z_{m+2k+1}$ for $m=0,1,\ldots, 2k$.
\end{prop}

\begin{proof}
The product in the lemma satisfies the equalities 
\begin{align*}
\prod_{i=1}^{2k-1}(1-u(z_{i-1})^*u(z_{i}))=\prod_{i=1}^{2k-1}\left(1+\frac{1}{2}(z_{i-1,+}+\bar{z}_{i-1,-})(z_{i,+}+\bar{z}_{i,-})\right)=\\
=\sum_{l=0}^{2k-1}\sum_{i_1<i_2<\ldots <i_l} 2^{-l}\prod_{j=1}^{l}\left((z_{i_j-1,+}+\bar{z}_{i_j-1,-})(z_{i_j,+}+\bar{z}_{i_j,-})\right).
\end{align*}
The Lemma follows from these equalities and degree reasons.
\end{proof}

Let us choose an open subset $U\subseteq Y$ such that there is a diffeomorphism $\nu:U\to B_{2n-1}$.  Let $\tilde{\nu}$ be as in equation \eqref{nuext} and define the function $\tilde{f}:\partial \Omega^{2k+1}\to \C$ by
\begin{equation}
\label{fsig}
\tilde{f}(z_0,z_1,\ldots, z_{2k}):= \sum_{\sigma\in S_{2(2k-1)}}\sum _{l=0}^{2k-1} (-1)^l2^{n-l-1}\epsilon_l(\sigma)\prod_{i=1}^l\langle \tilde{\nu} (f(z_{\sigma(2j-1)})),\tilde{\nu} (f(z_{\sigma(2j)}))\rangle
\end{equation}
where we identify $z_m$ with $z_{m+2k+1}$.

\begin{sats}
\label{holderdegree}
Suppose that $\Omega$ is a relatively compact strictly pseudo-convex domain with smooth boundary in a Stein manifold of complex dimension $n$ and that $Y$ is a connected, compact, orientable, Riemannian manifold of dimension $2n-1$. If $f:\partial \Omega\to Y$ is a H\"older continuous function of exponent $\alpha$ the degree of $f$ can be calculated by
\begin{align*}
\deg(f)=& (-1)^{n}\langle\tilde{\chi}_{\partial\Omega},g\circ f\rangle_k=\\
=&(-1)^{n}\int_{\partial\Omega^{2k+1}} \tilde{f} (z_0,z_1,\ldots, z_{2k})\prod _{j=0}^{2k}H_{\partial\Omega}(z_{j-1},z_j)\rd V=\\
=&(-1)^{n}\int_{\partial\Omega^{2k+1}} \tilde{f}(z_0,z_1,\ldots, z_{2k}) \prod _{j=0}^{2k}C_{\partial\Omega}(z_{j-1},z_j)\rd V
\end{align*}
whenever $2k+1>2n/\alpha$.
\end{sats}

\begin{proof}
By Theorem \ref{indexdegreeform} and Theorem \ref{holderindex} we have the equality
\[\deg(f)= (-1)^n\int_{\partial\Omega^{2k+1}} \tra\left(\prod_{j=0}^{2k}(1-g(f)(z_j)^*g(f)(z_{j+1}))H_{\partial\Omega}(z_{j-1},z_j)\right)\rd V.\]
Proposition \ref{nsch} implies 
\[\tra\left(\prod_{j=0}^{2k}(1-g(f)(z_j)^*g(f)(z_{j+1}))\right)=\tilde{f}(z_0,z_1,\ldots, z_{2k}),\]
from which the Theorem follows.
\end{proof}

Let us end this paper by a remark on the restriction in Theorem \ref{holderdegree} that the domain of $f$ must be the boundary of a strictly pseudo-convex domain in a Stein manifold. The condition on a manifold $M$ to be a a Stein manifold of complex dimension $n$ implies that $M$ has the same homotopy type as an $n$-dimensional $CW$-complex since the embedding theorem for Stein manifolds, see for instance \cite{grauert}, implies that a Stein manifold of complex dimension $n$ can be embedded in $\C^{2n+1}$ and by Theorem $7.2$ of \cite{milnor} an $n$-dimensional complex submanifold of complex euclidean space has the same homotopy type as a $CW$-complex of dimension $n$. 

Conversely, if $X$ is a real analytic manifold, then for any choice of metric on $X$, the co-sphere bundle $S^*X$ is diffeomorphic to the boundary of a strictly pseudo-convex domain in a Stein manifold, see for instance Proposition $4.3$ of \cite{guille} or Chapter V.$5$ of \cite{grauert}. So the degree of $f$ coincides with the mapping $H^{2n-1}_{dR}(S^*Y)\to H^{2n-1}_{dR}(S^*X)$ that $f$ induces under the Thom isomorphism $H^n_{dR}(X)\cong H^{2n-1}(S^*X)$. Thus the degree of a function $f:X\to Y$ can be expressed using our methods for any real analytic $X$.

\end{document}